\magnification 1200
\input amssym.def
\input amssym.tex
\parindent = 40 pt
\parskip = 12 pt

\font \heading = cmbx10 at 12 true pt
 at 22 true pt
\font \medheading =cmbx7 at 16 true pt
 at 7 true pt

\def \R{{\bf R}}

\centerline{\medheading Fourier transforms of powers of well-behaved}
\centerline{\medheading 2D real analytic functions}
\rm
\line{}
\line{}
\centerline{\heading Michael Greenblatt}
\line{}
\centerline{May 25, 2016}
\baselineskip = 12 pt
\font \heading = cmbx10 at 15 true pt
 at 13 true pt
\line{}
\line{}
\line{} 
\noindent{\heading 1. Introduction and Definitions.}

In this paper we consider Fourier transforms of powers of local two-dimensional real analytic functions. Namely we
consider integrals
$$F(\lambda_1,\lambda_2) = \int_{\R^2} \phi(x_1,x_2)|f(x_1,x_2)|^{-\rho}e^{-i\lambda_1x_1 - i\lambda_2x_2}\,dx_1\,dx_2 \eqno (1.1)$$
Here $f(x_1,x_2)$ is real analytic near the origin with $f(0,0) = 0$, $\rho >0$ such that $|f(x_1,x_2)|^{-\rho}$ is integrable on a
neighborhood of the origin, and $\phi(x_1,x_2)$ is supported on a neighborhood of the origin and  $C^1$ on $(\R - \{0\})^2$ such
 that for some constant $A > 0$, on $(\R - \{0\})^2$  we have
$$|\phi(x_1,x_2)| < A {\hskip 0.8 in} |\nabla \phi(x_1,x_2)| < {A \over (x_1^2 + x_2^2)^{1 \over 2}}    \eqno (1.2)$$
 The prototype for $\phi(x_1,x_2)$ would be a smooth cutoff 
function on a neighborhood of the origin, but since the more general form is no more difficult we stipulate this condition. Note that
we can multiply $\phi(x_1,x_2)$ by the characteristic function of any quadrant and $(1.2)$ still holds. This allows us for example
to estimate the Fourier transform of $|f(|x_1|, |x_2|)|^{-\rho}$ by adding the estimates for a given quadrant.

In the paper [G4], we provide various sharp estimates that can be proven for the functions $F(\lambda_1,\lambda_2)$.
The theorems of [G4] are stated in a rather general form, and as a result sometimes the estimates of that paper are 
not amenable to being written out directly in terms of explicit properties of $f(x_1,x_2)$. 

In this paper, we will expand on the results of [G4] and define a class of "well-behaved" functions that contains a
 number of relevant examples for which such estimates can be explicitly described. Specifically, we
will see that for a range of $\rho$, for these well-behaved $f(x_1,x_2)$ we will be able to find optimal estimates of the form $|F(\lambda)| < C|\lambda_i|^{-\epsilon_i}$ for $i = 1, 2$, which immediately lead to optimal estimates of the form 
$|F(\lambda)| < C|\lambda|^{-\epsilon}$. Here
$\lambda$ denotes $(\lambda_1,\lambda_2)$. The $\epsilon_i$ will be explicitly describable in terms of the Newton polygons of
$f(x_1,x_2)$.  We will further see that for a subclass of these functions, these estimates hold
for all $\rho$ (even when $\rho < 0$) and furthermore we even have
 estimates $|F(\lambda)| < C\alpha(\lambda)$, where again the estimates can be explicitly expressed in
terms of the Newton polygon of $f(x_1,x_2)$. 

In order to state our theorems, we now give some terminology that is frequently used in the subject of two-dimensional oscillatory
integrals. 

\noindent {\bf Definition 1.1.} Let $f(x_1,x_2) = \sum_{a,b} f_{ab}x^ay^b$ denote the 
Taylor expansion of $f(x_1,x_2)$ at the origin. For any $(a,b)$ for which $f_{ab} \neq 0$, let $Q_{ab}$ be the
quadrant $\{(x_1,x_2) \in \R^2: 
x \geq a, y \geq b \}$. Then the {\it Newton polygon} $N(f)$ of $f(x_1,x_2)$ is defined to be 
the convex hull of the union of all $Q_{ab}$.  

In general, the boundary of $N(f)$ consists of finitely many (possibly none) bounded edges of negative slope
as well as an unbounded vertical ray and an unbounded horizontal ray. 

\noindent A key role in our paper is played by the following polynomials.

\noindent {\bf Definition 1.2}. Suppose $e$ is a compact edge of $N(f)$. Define $f_e(x_1,x_2)$
by $f_e(x_1,x_2) = \sum_{(a,b) \in e} f_{ab} x^{a}y^{b}$. In other words $f_e(x_1,x_2)$ is the sum of the terms
of the Taylor expansion of $f$ corresponding to $(a,b)$ on the edge $e$.

\noindent It is an important point that one only considers compact edges of $N(f)$ in the above definition. Next, our theorem statements will make use the following notion.

\noindent {\bf Definition 1.3.} The {\it Newton distance} $d(f)$ of $f(x_1,x_2)$ is defined to be 
$\inf \{t: (t,t) \in N(f)\}$.

\noindent Our well-behavedness condition is then given by the following.

\noindent {\bf Definition 1.4.} $f(x_1,x_2)$ is said to be {\it well-behaved} if the order of any zero of any $f_e(x_1,x_2)$ 
in $(\R - \{0\})^2$ is less than $d(f)$, and if there is an edge $e$ of slope $-1$ then that $f_e(x_1,x_2)$  has no zeroes at all in
$(\R - \{0\})^2$.

This condition is related to the concept of adapted coordinates in the subject of two-dimensional oscillatory integrals as 
 initiated in [V]. Namely, $f(x_1,x_2)$ is in adapted coordinates if the zeroes of each $f_e(x_1,x_2)$ in $(\R - \{0\})^2$ have
order less than or equal to $d(f)$. It turns out that that the scalar oscillatory index of $f(x_1,x_2)$ at the origin (see [AGV] for
the relevant definitions) is equal to ${1 \over d(f)}$ if and only if $f(x_1,x_2)$ is in adapted coordinates, and thus
 in this situation one can readily compute this index in terms of $N(f)$. The reference  [AGV] has a wealth of information
 on related matters. For the purpose of this paper, we are most concerned with the following (closely related) fact.

\noindent {\bf Lemma 1.1. ([G1])} If $f(x_1,x_2)$ is well-behaved, then $|f(x_1,x_2)|^{-\rho}$ is integrable on a neighborhood of the origin
whenever $\rho < {1 \over d(f)}$, and is not integrable on any neighborhood of the origin whenever $\rho > {1 \over d(f)}$. 

 Next, we define the function $f^*(x_1,x_2)$, which will be a regularized version of $f(x_1,x_2)$ whose general behavior will
be the same as $f(x_1,x_2)$ when $f(x_1,x_2)$ is well-behaved but for which many relevant quantities such as integrals are quite a bit easier to compute.

\noindent {\bf Definition 1.5.} $f^*(x_1,x_2)$ denotes the function $\sum_{(v_1,v_2)\,\,a\,\, vertex\,\, of N(f)} |x_1|^{v_1}|x_2|^{v_2}$.

\noindent A useful fact concerning $f^*(x_1,x_2)$ is the following.

\noindent {\bf Lemma 1.2.} Suppose $0 < \rho < {1 \over d(f)}$ and $f(x_1,x_2)$ is well-behaved.
 Then there are positive constants $C_1$ and $C_2$ depending on $\rho$ and $f$ such that  if $R$ is a dyadic rectangle one has
$$C_1 \int_R |f|^{-\rho} \leq \int_R (f^*)^{-\rho}  \leq C_2  \int_R |f|^{-\rho} \eqno (1.3)$$

\noindent {\bf Proof.} The $n$-dimensional version of this was proven in [G2]. Specifically, by Lemma 2.1 of [G2], one has
the existence of a constant $C$ for which $|f(x)| \leq Cf^*(x)$ for all $x$, which gives the right-hand side of $(1.3)$. The 
left-hand side follows from $(4.15)$ of [G3], taking $\epsilon = 1$, since the left hand inequality of $(1.3)$ for the portion of $R$
where $|f(x)| > f^*(x)$ is immediate. 

As a consequence of  Lemmas 1.1-1.2, if $f(x_1,x_2)$ is well-behaved $(f^*(x_1,x_2))^{-\rho}$ is integrable on a neighborhood of the origin if $\rho  < {1 \over d(f)}$
and integrable on no neighborhood of the origin if $\rho  > {1 \over d(f)}$ since the same is true for $|f(x_1,x_2)|$. It can be shown
that the same is also true for $f(x_1,x_2)$ that is not well-behaved.

\noindent{\heading 2. Main Results.}

 Our first lemma defines some quantities used in the 
statement of our first theorem. In the following, we denote the edges of $N(f)$  by $e_0,...,e_n$,
where $e_0$ is the horizontal edge, $e_n$ is the vertical edge, and the $e_i$ are listed in order of decreasing slope. We write
the slope of $e_i$ as $-{1 \over m_i}$ where for $i = 0$ we take $m_i = \infty$ and for $i = n$ we take 
$m_i = 0$. Thus $m_{i+1} < m_i$ for all $i$. We denote by $v_i$ the vertex of $N(f)$ between edges $e_{i-1}$ and $e_i$, and write $v_i = (v_1^i,v_2^i).$

\noindent {\bf Lemma 2.1.} Suppose that $(f^*(x_1,x_2))^{-\rho}$ is integrable on a neighborhood of the 
origin, where $\rho > 0$. Then there exists an $\epsilon \geq 0$ and a $d = 0$ or $1$, both depending on $\rho$ and $f$, 
 such that if $0 <  r_0 <  {1 \over 2}$ there are positive constants $c_1, c_2$ depending on $\rho, f$, and $r_0$, such that
if $0 < r < r_0$ one has
one has $c_1  r^{-\epsilon}|\ln r|^d \leq \int_0^{r_0}  (f^*(r,x_2))^{-\rho}\,dx_2 \leq c_2  r^{-\epsilon}|\ln r|^d$.

\noindent In addition, $\epsilon$ and $d$ can be explicitly computed by finding the dominant term in the sum
$$r^{-\rho v_1^1}\int_0^{r^{m_1}}x_2^{-\rho v_2^1}\,dx_2 + \sum_{i=2}^{n-1} r^{-\rho v_1^i}\int_{r^{m_{i-1}}}^{r^{m_i}} x_2^{-\rho v_2^i}\,dx_2  +  r^{-\rho v_1^n}\int_{r^{m_{n-1}}}^1  x_2^{-\rho v_2^n}\,dx_2\eqno (2.1)$$
In the event that $n = 2$, one excludes the middle term of $(2.1)$, and in the event $n = 1$ we replace $(2.1)$ by
$r^{-\rho v_1^1}\int_0^1 x_2^{-\rho v_2^1}\,dx_2$.

\noindent {\bf Proof.} We first consider the portion of the integral $\int_0^r (f^*(r,x_2))^{-\rho}\,dx_2 $
between any $x_2 = r^{m_{i-1}}$ and $x_2 = r^{m_i}$ for  $i \geq 2$ that occurs. We claim that in this range, the quantity $r^{{v_1}^i}
x_2^{{v_2}^i}$ is at least as large as $r^{{v_1}^j} x_2^{{v_2}^j}$ for any $j \neq i$. To see why this is the case, we
look at the ratio $(r^{{v_1}^i}x_2^{{v_2}^i})/(r^{{v_1}^j} x_2^{{v_2}^j}) = r^{v_1^i -v_1^j}x_2^{v_2^i - v_2^j}$.
If $v_2^i >  v_2^j$, since $x_2 \geq  r^{m_{i-1}}$  we have 
$$r^{v_1^i -v_1^j}x_2^{v_2^i - v_2^j} \geq   r^{v_1^i -v_1^j +m_{i-1}(v_2^i - v_2^j)} \eqno (2.2)$$
Because  $(v_1^j, v_2^j)$ is on or above the edge $e_{i-1}$, whose slope is $-{1 \over m_{i-1}}$, we have that $v_2^j + m_{i-1} v_1^j \geq v_2^i + m_{i-1}v_1^i$. Thus the exponent in $(2.2)$ is negative. Hence $r^{{v_1}^i}
x_2^{{v_2}^i} \geq  r^{{v_1}^j} x_2^{{v_2}^j}$ as needed. If on the other hand $v_2^i \leq v_2^j$, since 
$x_2  \leq r^{m_i}$, in place of $(2.2)$ we can use
$$r^{v_1^i -v_1^j}x_2^{v_2^i - v_2^j}  \geq  r^{v_1^i -v_1^j +m_i(v_2^i - v_2^j)} \eqno (2.3)$$
This time, we use that since $(v_1^j, v_2^j)$ is on or above the edge $e_i$ we have $v_2^j + m_iv_1^j \geq v_2^i + m_iv_1^i$. Thus the exponent in $(2.3)$ is again negative and  $r^{{v_1}^i} x_2^{{v_2}^i} \geq  r^{{v_1}^j} x_2^{{v_2}^j}$ as desired.

Hence we have seen that on the portion of the integral where $x_1^{m_{i-1}} \leq x_2 \leq x_1^{m_i}, i \geq 2$, we have that
$r^{{v_1}^i}x_2^{{v_2}^i} \geq r^{{v_1}^j} x_2^{{v_2}^j}$ for $j \neq i$. Thus $f^*(r,x_2)$ is the sum of several positive terms, the largest of which is $r^{{v_1}^i}x_2^{{v_2}^i}$. Hence there are constants $C_1$ and $C_2$ depending on 
$N(f)$ and $\rho$ such that whenever $x_1^{m_i} \leq  x_2  \leq x_1^{m_{i-1}}$ we have
$$C_1 r^{{v_1}^i}x_2^{{v_2}^i} < f^*(r,x_2) < C_2 r^{{v_1}^i}x_2^{{v_2}^i} \eqno (2.4)$$
Next, we will prove an analogue of $(2.4)$ that holds on $x_2 < x_1^{m_1}$. This time $i = 1$ is the dominant term. To see why, 
note that since $v_2^1 \leq v_2^j$ for all $j \neq 1$ and $x_2 < x_1^{m_1}$, we have that $(2.3)$ holds for $i = 1$ 
and all $j \neq 1$. So since each $(v_1^j, v_2^j)$ is on or above the edge $e_1$, we have $v_2^j + m_1v_1^j \geq v_2^1 + m_1v_1^1$ and like before we have $r^{{v_1}^1} x_2^{{v_2}^1} \geq  r^{{v_1}^j} x_2^{{v_2}^j}$. The analogue to $(2.4)$
that we get for the points were $y < x^{m_1}$ is therefore
$$C_1' r^{{v_1}^1}x_2^{{v_2}^1} < f^*(r,x_2) < C_2'  r^{{v_1}^1}x_2^{{v_2}^1} \eqno (2.5)$$
Similarly, if $x_2 > r^{m_{n-1}}, n \geq 2$, one can argue as in the above cases and show that we have
$$C_1'' r^{{v_1}^n}x_2^{{v_2}^n} < f^*(r,x_2) < C_2''  r^{{v_1}^n}x_2^{{v_2}^n} \eqno (2.6)$$
Equations $(2.4)-(2.6)$  cover the entire $y$ range of integration in  $\int_0^{r_0} (f^*(r,x_2))^{-\rho}\,dx_2 $, except when $N(f)$ has
exactly one vertex. But in this case $f^*(r,x_2) =  r^{{v_1}^1}x_2^{{v_2}^1}$ which serves as a substitute for $(2.4)-(2.6)$.

Equation $(2.1)$ follows from $(2.4)-(2.6)$ in short order; one simply takes the monomial from $(2.4)-(2.6)$, raises it to the
$-\rho$ power, and integrates in $x_2$ over its domain. Adding over all domains gives $(2.1)$. This completes the proof of 
Lemma 2.1.

Note that by the proof of Lemma 2.1, one has that $f^*(x_1,x_2)$ is always within a constant factor of some
 dominant $x_1^{v_1^i}x_2^{v_2^i}$ which can be readily determined at a given $(x_1,x_2)$. This will prove useful later.

Note also that the $\epsilon$ given by the expression $(2.1)$ is a continuous function of $\rho$ at any value of $\rho$ where the
expression is finite. As a result, when $(2.1)$ is finite for $\rho = {1 \over d(f)}$ this
$\epsilon$ must be $1$. This true for the following reason.  Since $|f^*(x_1,x_2)|^{-\rho}$ is integrable  on a neighborhood of the origin 
when $\rho < {1 \over d(f)}$ by Lemma 1.1,  $\epsilon$ must be less than $1$ for such $\rho$. By continuity
 $\epsilon$ is therefore at most $1$ when $\rho ={1 \over d(f)}$.
If it were strictly less, the continuity of $\epsilon$ in the expression $(2.1)$ implies that we could integrate  $|f^*(x_1,x_2)|^{-\rho}$ to a finite value on a neighborhood of the origin for some $\rho > {1 \over d(f)}$, which is not possible by Lemma 1.1.
Hence $\epsilon = 1$ when $\rho = {1 \over d(f)}$ whenever $(2.1)$ is finite.

As a result, the continuity of $\epsilon$ in $\rho$ says there will be an interval
on which $\epsilon > {1 \over 2}$ as long as $(2.1)$ is finite at $\rho = {1 \over d(f)}$, which is the typical situation (but not always;
see Example 1 below.) By symmetry the same will be true when $f^*(x_1,x_2)$ is replaced by $f^*(x_2,x_1)$.
 This justifies the $\epsilon_i > {1 \over 2}$ conditions in the statement of Theorem 2.2, our first main theorem, which we now
come to.

\noindent {\bf Theorem 2.2.} Suppose $\rho > 0$ and 
$f(x_1,x_2)$ is well-behaved and the $\epsilon$ of Lemma 2.1 is greater than ${1 \over 2}$ for both $f^*(x_1,x_2)$ and $f^*(x_2,x_1)$.

\noindent {\bf a)} Let $(\epsilon_1, d_1)$ be as in Lemma 2.1 applied to $f(x_1,x_2)$, and let $(\epsilon_2, d_2)$ as in Lemma 2.1 applied to $f(x_2,x_1)$. Then there is a neighborhood $N$ of the origin such that if the function $\phi(x_1,x_2)$ in $(1.1)$ is supported in $N$ then we have the following
estimates, where $C$ is a constant depending on $f$, $\rho$, $N$, and the constant $A$ of $(1.2)$.
$$|F(\lambda_1,\lambda_2)| < C(2 + |\lambda_1|)^{\epsilon_1 - 1}(\ln (2 + |\lambda_1|))^{d_1} \eqno (2.7a)$$
$$ |F(\lambda_1,\lambda_2)| < C(2 + |\lambda_2|)^{\epsilon_2 - 1}(\ln (2 + |\lambda_2|))^{d_2} \eqno (2.7b)$$
Thus $ |F(\lambda_1,\lambda_2)| < C(2 + |\lambda|)^{\epsilon_j - 1}(\ln (2 + |\lambda|))^{d_j}$ where $(\epsilon_j, d_j)$ 
denotes the slower of the two decay rates. 

\noindent {\bf b)} When $\phi(x_1,x_2)$ is bounded below by a positive constant on a neighborhood  of the origin, then the 
exponents $\epsilon_i$ of $(2.7a)-(2.7b)$ are best possible whenever $\epsilon_i < 1$; one does 
not have an estimate $|F(\lambda_1,\lambda_2)| < C(2 + |\lambda_i|)^{\epsilon' - 1}$ for $\epsilon'< \epsilon_i$.

Although we won't prove it here, there is a variation of Theorem 2.2 for the case when $f(x_1,x_2)$ is not well-behaved, but 
where instead the sum of terms of lowest degree  of $f(x_1,x_2)$ has no zeroes on $V = \{(x_1,x_2): |x_2| < c|x_1|\} \cap (\R - \{0\})^2$ for some 
$c > 0$. If instead of $F(\lambda_1,\lambda_2)$ one looks at the Fourier transform of $|f(x_1,x_2)|^{-\rho} \chi_V(x_1,x_2)$, one can
show that Theorem 2.2 holds where the new $(\epsilon_1,d_1)$ and $(\epsilon_2,d_2)$ are defined by the variation on Lemma 2.1
 where the integral $\int_0^{r_0} (f^*(r,x_2))^{-\rho}\,dx_2$ in the statement  is replaced by $\int_0^{c x_1} (f^*(r,x_2))^{-\rho}\,dx_2$ and 
 where the integral $\int_0^1 (f^*(x_1,r))^{-\rho}\,dx_2$ in the statement  is replaced by $\int_{x_2 \over c}^1 (f^*(x_1,r))^{-\rho}\,dx_1$. This variant of Theorem 2.2 allows us to divide a neighborhood
of the origin via lines through the origin, resulting in several wedges $W_i$. One can rotate each $W_i$ to turn it into a set of
the form $V$. If this variation of Theorem 2.2 applies on each such $V$, then can estimate $F(\lambda_1,\lambda_2)$ by adding 
the Fourier transform estimates for each  $|f(x_1,x_2)|^{-\rho} \chi_V(x_1,x_2)$.

\noindent {\bf Example 1.} Suppose $f(x_1,x_2) = x_1^ax_2^b$ for some $a$ and $b$ not both zero. Then 
$|f(x_1,x_2)|^{-\rho}$ is
integrable on a neighborhood of the origin if $\rho < {1 \over max(a,b)}$. Here $f^*(x_1,x_2) = |x_1|^a|x_2|^b$, and one can
compute $(\epsilon_1,d_1)$ using the integral $\int_0^1 (r^ax_2^b)^{-\rho}\,dx_2  = Cr^{-a\rho}$. So $\epsilon_1 = a\rho$ 
and $d_1 = 0$ here. By symmetry, $\epsilon_2 = b\rho$ and $d_2 = 0$. Hence Theorem 2.2 says that if 
${1\over 2\min (a,b)} < \rho < {1 \over \max(a,b)}$
one has estimates $|F(\lambda_1,\lambda_2)| \leq C|\lambda_1|^{a\rho - 1}$ and $|F(\lambda_1,\lambda_2)| \leq C|\lambda_2|^{b\rho - 1}$, leading to an overall decay rate of $|F(\lambda_1,\lambda_2)| \leq C|\lambda|^{\max(a\rho, b\rho) - 1}$.
Note that if $\min(a,b)  \leq {1 \over 2}\max(a,b)$, the conditions of Theorem 2.2 will never hold for this example.

\noindent {\bf Example 2.} Suppose $f(x_1,x_2) = |x_1|^a + |x_2|^b$ for some $a$ and $b$ neither of which is zero. Then 
$f^*(x_1,x_2) = |x_1|^a + |x_2|^b$. The Newton polygon $N(f)$ has two vertices, $(a,0)$ and $(0,b)$, and three edges:
the vertical and horizontal edges, and a compact edge of slope $-{b \over a}$. So $m_0 = \infty$, $m_1 = {a \over b}$, and $m_2 = 0$.  Then $(\epsilon_1,d_1)$ is computed using
$$r^{-\rho a}\int_0^{r^{a \over b}}1 \,dx_2 + \int_{r^{a \over b}}^1 \,\,x_2^{-\rho b}\,dx_2$$
If $\rho \neq {1 \over b}$, this is equal to $r^{-\rho a + {a \over b}} + {1 \over 1 - \rho b}(1 - r^{(1 - \rho b) {a \over b}}) =
{1 \over 1 - \rho b} -{\rho b \over 1 - \rho b} r^{-\rho a + {a \over b}}$. The second term dominates if $\rho  > {1 \over b}$, 
and the first term dominates if $\rho < {1 \over b}$. Hence $(\epsilon_1, d_1) = (\rho a - {a \over b}, 0)$ if $\rho > {1 \over b}$ and
$(\epsilon_1, d_1) = (0,0)$ if $\rho < {1 \over b}$. If $\rho = {1 \over b}$, the sum of the two integrals is 
$r^{-\rho a + {a \over b}} - {a \over b} \ln r$. Thus the second term dominates, and $(\epsilon_1,d_1) = (0, 1)$.
By symmetry, $(\epsilon_2, d_2) = (\rho b - {b \over a}, 0)$ if $\rho > {1 \over a}$, $(\epsilon_2, d_2) = (0,0)$ if 
$\rho < {1 \over a}$, and $(\epsilon_2, d_2) = (0,1)$ if $\rho = {1 \over a}$.

$f(x_1,x_2)$ is integrable over a neighborhood of the origin if $\epsilon_1 < 1$, which in the current situation is  equivalent to the statement that
$\epsilon_2 < 1$. The condition works out to $\rho < {1 \over a} + {1 \over b}$. One has $\epsilon_1 > {1 \over 2}$ if 
$\rho a - {a \over b} > {1 \over 2}$ or equivalently $\rho > {1 \over 2a} + {1 \over b}$. Similarly, one has $\epsilon_2 
> {1 \over 2}$ if $\rho > {1 \over a} + {1 \over 2b}$. Hence Theorem 2.2 applies for $\rho$ on the smaller of the two intervals
 $({1 \over 2a} + {1 \over b}, {1 \over a} + {1 \over b})$ or $({1 \over a} + {1 \over 2b}, {1 \over a} + {1 \over b})$. For
such $\rho$ one has estimates $|F(\lambda_1,\lambda_2)| \leq C|\lambda_1|^{\rho a - {a \over b} - 1 }$ and $|F(\lambda_1,\lambda_2)| \leq C|\lambda_2|^{\rho b - {b \over a} - 1}$. The overall estimate obtained is then $|F(\lambda_1,\lambda_2)| 
\leq C|\lambda|^{\max( \rho a - {a \over b} - 1, \rho b  - {b \over a} - 1) }$. If one works it out, one sees that one uses the exponent
$ \rho a - {a \over b} - 1 $ if $a \leq b$ and the exponent $\rho b  - {b \over a} - 1$ if $a \geq b$. 

This example also satisfies the conditions of Theorem 2.3. As a result the estimate  $|F(\lambda_1,\lambda_2)| 
\leq C|\lambda|^{\max( \rho a - {a \over b} - 1, \rho b  - {b \over a} - 1) }$ (as well as more precise estimates) will hold for any $(a,b)$ and any $\rho$.

Our second main theorem will give more precise information than Theorem 2.2 when each $f_e(x_1,x_2)$ has no zeroes on
$(\R - \{0\})^2$. Instead of $(1.2)$ we will assume that $\phi(x_1,x_2)$ is $C^{\infty}$ on 
$(\R - \{0\})^2$ and there are constants $A$ and $A_{a,b}$ such that 
$$|\phi(x_1,x_2)| \leq A {\hskip 0.8 in}|\partial_{x_1}^a\partial_{x_2}^b \phi(x_1,x_2)| \leq A_{a,b}|x_1|^{-a}|x_2|^{-b} {\hskip 0.2 in} (\forall a \,\forall b) \eqno (2.8)$$
\noindent {\bf Theorem 2.3.} Suppose $(2.8)$ holds and that each $f_e(x_1,x_2)$ 
has no zeroes on $(\R - \{0\})^2$. Then there is a neighborhood $U$ of the origin such that if $\phi(x_1,x_2)$ is supported in $U$, then equations  $(2.7a)-(2.7b)$ hold for all $\rho$, even if $\rho < 0$. In fact, one has the following stronger estimate, where 
$C$ is a constant depending on $f$, $\rho$, $U$, and the constant $A$ of $(1.2)$.
$$|F(\lambda_1,\lambda_2)| \leq C\int_{\{(x_1,x_2) \in U:\,|x_1| < |\lambda_1|^{-1}, \,|x_2| < |\lambda_2|^{-1}\}}|f^*(x_1,x_2)|^{-\rho}\,dx_1\,dx_2  \eqno (2.9)$$

Up to a constant factor, one can explicitly determine the integral $(2.9)$ similarly to in Theorem 2.2, proceeding as in Lemma 2.1
where one divides a neighborhood of the origin into domains on each of which 
$f^*(x_1,x_2)$ is within a constant factor of some explicitly determinable $x_1^{v_1^i}x_2^{v_2^i}$. On each such domain,
the right-hand side of $(2.9)$ will be within a bounded factor of $|\lambda_1|^a|\ln \lambda_1|^{d_1} |\lambda_2|^b|\ln \lambda_2|^{d_2}$
 for some $a$ and $b$ and $d_i = 0$ or $1$. 

It can also be shown similarly to the proof of Theorem 2.2b) 
that if $a < 0$ and $b < 0$ then on the domain
 $\{(\lambda_1,\lambda_2): |\lambda_1|^{-m_i} < |\lambda_2|^{-1} <|\lambda_1|^{-m_{i+1}} \}$
 for compact edges $e_i$ and $e_{i+1}$, the exponents $a$ and $b$ are best possible.

\noindent {\bf Example 1.} Let $f(x_1,x_2) = x_1^ax_2^b$ where $a$ and $b$ are not both zero. Then the right-hand side of $(2.9)$ 
is given by $C\int_0^{|\lambda_1|^{-1}}\int_0^{|\lambda_2|^{-1}} x_1^{-a\rho}x_2^{-b\rho}\,dx_2\,dx_1$, which equals
 $C|\lambda_1|^{a\rho - 1}|\lambda_2|^{b\rho - 1}$ when it is finite. Note the improvement over the estimate for the same example after 
Theorem 2.2.

\noindent {\bf Example 2.} Suppose $f(x_1,x_2) = |x_1|^a + |x_2|^b$ for some $a$ and $b$ neither equal to zero. 
Theorem 2.3 then gives
$$|F(\lambda_1,\lambda_2)| \leq C\int_0^{|\lambda_1|^{-1}}\int_0^{|\lambda_2|^{-1}}( x_1^a + x_2^b)^{-\rho}\,dx_2\,dx_1 \eqno (2.10)$$
One divides the integral along the curve $x_2 = x_1^{a \over b}$, and $(2.10)$ becomes
$$|F(\lambda_1,\lambda_2)| \leq C\int_0^{|\lambda_1|^{-1}}\int_0^{\min(x_1^{a \over b}, |\lambda_2|^{-1})}x_1^{-a\rho}
\,dx_2\,dx_1$$
$$ + C\int_0^{|\lambda_1|^{-1}}\int_{\min(x_1^{a \over b}, |\lambda_2|^{-1})}^{ |\lambda_2|^{-1}} x_2^{-b\rho}\,dx_2 \,dx_1\eqno (2.11)$$
One can readily perform the integrations in $(2.11)$ to get explicit formulas. One gets two different formulas depending on
whether or not $|\lambda_1|^{-a} \leq |\lambda_2|^{-b}$.

\noindent {\heading 3. Theorem proofs.}

\noindent {\bf Proof of Theorem 2.2.} 

\noindent Suppose we are in the setting of Theorem 2.2. The key fact that we use here is Corollary 3.4 of [G4], which says 
that 
$$|F(\lambda_1,\lambda_2)|  \leq C \int_N (1 + |\lambda_1 x_1| + |\lambda_2 x_2|)^{-{1 \over 2}} |f(x_1,x_2)|^{-\rho}\,dx_1\,dx_2 \eqno (3.1)$$
Here $N$ is a small neighborhood of the origin on which the resolution of singularities algorithm of [G5] applies, and we
henceforth assume $\phi(x_1,x_2)$ is supported on $N$. Thus for $i = 1,2$ we have
$$|F(\lambda_1,\lambda_2)|  \leq C \int_N (1 + |\lambda_i x_i|)^{-{1 \over 2}} |f(x_1,x_2)|^{-\rho}\,dx_1\,dx_2 \eqno (3.2)$$
Suppose $|\lambda_i| > 2$. Splitting the integral $(3.2)$ at $|x_i| = {1 \over |\lambda_i|}$, equation $(3.2)$ becomes
$$|F(\lambda_1,\lambda_2)|  \leq C \int_{\{(x_1,x_2) \in N: |x_i| < {1 \over |\lambda_i|}\}} |f(x_1,x_2)|^{-\rho}\,dx_1\,dx_2 $$
$$+ C{1 \over |\lambda_i|^{1 \over 2}} \int_{\{(x_1,x_2) \in N: |x_i| \geq  {1 \over |\lambda_i|}\}} {x_i}^{-{1 \over 2}}    |f(x_1,x_2)|^{-\rho}\,dx_1\,dx_2       \eqno (3.3)$$
By Lemma 1.2, one can replace $|f(x_1,x_2)|$ by $f^*(x_1,x_2)$ in $(3.3)$. In the two integrals of the resulting expression, we
first integrate in the variable that is not $x_i$, inserting the right-hand inequality of Lemma 2.1. The result is
$$|F(\lambda_1,\lambda_2)|  \leq C \int_0^{1 \over |\lambda_i|} {x_i}^{-\epsilon_i}|\ln x_i|^{d_i}\,dx_i + 
C{1 \over |\lambda_i|^{1 \over 2}}\int_{1 \over |\lambda_i|}^{1 \over 2} x_i^{-{1 \over 2} - \epsilon_i}|\ln x_i|^{d_i}\,dx_i$$
Integrating the two terms in $(3.3)$ and using that $\epsilon_i > {1 \over 2}$ we obtain the desired estimate
$$|F(\lambda_1,\lambda_2)| \leq C |\lambda_i|^{\epsilon_i - 1}|\ln \lambda_i|^{d_i} \eqno (3.4)$$
This is $(2.7a)-(2.7b)$ when $|\lambda_i| > 2$. When $|\lambda_i| < 2$, one obtains $(2.7a)-(2.7b)$ simply by taking absolute
values and integrating. Thus we have proved $(2.7a)-(2.7b)$ and the proof of part a) of Theorem 2.2 is complete. 

Moving on to part b), suppose $\epsilon_i > 0$, $\phi(x_1,x_2)$ is bounded below on a neighborhood of the 
origin and the estimate $|F(\lambda_1,\lambda_2)| \leq C(2 + |\lambda_i|)^{\epsilon' - 1}$ holds, where $\epsilon' < \epsilon_i < 1$,
and we will reach a contradiction. Without
loss of generality we take $i = 1$. Let $\psi(x)$ be a smooth function on 
$\R$ whose Fourier transform is a nonnegative compactly supported function equal to $1$ on a neighborhood of the origin. Since $\epsilon' < \epsilon_1< 1$
and $\epsilon_1 > 0$, we may let $\eta > 0$ be 
 such that $0 < \eta + \epsilon' < \epsilon_1$. For a large $K$ we look at
$$I_K = \int F(\lambda_1,0) \psi\bigg({\lambda_1 \over K}\bigg)|\lambda_1|^{-\eta -\epsilon'}\,d\lambda_1 \eqno (3.5)$$
Since  $|F(\lambda_1,0)| \leq C(2 + |\lambda_1|)^{\epsilon' - 1}$, we have that
$$|I_K| \leq C\int (2 + |\lambda_1|)^{\epsilon' - 1}|\lambda_1|^{-\eta -\epsilon'}\,d\lambda_1 \eqno (3.6)$$
Because $\eta > 0$, the integrand in $(3.6)$  is integrable for large $|\lambda_1|$, and because $\eta  + \epsilon' < 1$ the 
integrand in $(3.6)$ is integrable for small $|\lambda_1|$. Hence the $I_K$ are uniformly bounded in $K$. On the other hand
$$I_K = \int_{\R^2} \phi(x_1,x_2)|f(x_1,x_2)|^{-\rho}e^{-i\lambda_1x_1} \psi\bigg({\lambda_1 \over K}\bigg)|\lambda_1|^{-\eta -\epsilon'}\,d\lambda_1\,dx_1\,dx_2 \eqno (3.7)$$
Performing the $\lambda_1$ integral in $(3.7)$ leads to
$$I_K = \int_{\R^2} \phi(x_1,x_2)|f(x_1,x_2)|^{-\rho} K^{1 -\eta - \epsilon '} \xi(K x_1)\,dx_1\,dx_2 \eqno (3.8)$$
Here $\xi$ is the Fourier transform of $\psi(\lambda_1)|\lambda_1|^{-\eta - \epsilon '}$. Since the Fourier transform of 
$\psi(\lambda_1)$ is nonnegative and the Fourier transform of $|\lambda_1|^{-\eta - \epsilon '}$ is of the form
 $c|x_1|^{\eta + \epsilon ' - 1}$, $\xi(x_1)$ is of the form $c \tilde{\xi}(x_1)$ where  $\tilde{\xi}(x_1)$ is nonnegative and
decays as $|x_1|^{\eta + \epsilon ' - 1}$ as $|x_1| \rightarrow \infty$. Thus we can rewrite $(3.8)$ as 
$$|I_K| = |c| \int_{\R^2} \phi(x_1,x_2)|f(x_1,x_2)|^{-\rho} K^{1 -\eta - \epsilon '} \tilde{\xi}(K x_1)\,dx_1\,dx_2 \eqno (3.8')$$
Since $\phi(x_1,x_2)$ is nonnegative
and is positive on a neighborhood of the origin, there is a constant $C$ and a neighborhood $N$ of the origin such that
$$I_K \geq C\int_{N} |f(x_1,x_2)|^{-\rho} K^{1-\eta - \epsilon '} \tilde{\xi}(K x_1)\,dx_1\,dx_2 \eqno (3.9)$$
Shrinking $N$ if necessary and assuming $N$ is a union of dyadic rectangles on which Lemma 1.2 holds, we therefore have
$$I_K \geq C'\int_{N} |f^*(x_1,x_2)|^{-\rho} K^{1 -\eta - \epsilon '} \tilde{\xi}(K x_1)\,dx_1\,dx_2 \eqno (3.10)$$
Performing the $x_2$ integration and using Lemma 2.1, for some $a > 0$ we therefore have
$$I_K \geq C''\int_{-a}^a K^{1 -\eta - \epsilon '}\tilde{\xi}(K x_1) x_1^{-\epsilon_1}|\ln x_1|^{d_1}  \,dx_1 \eqno (3.11)$$
In particular, for any $b > 0$ we have
$$I_K \geq  C''\int_{b}^a K^{1 -\eta - \epsilon '}\tilde{\xi}(K x_1) x_1^{-\epsilon_1}|\ln x_1|^{d_1}  \,dx_1 \eqno (3.12)$$
Taking limits as $K \rightarrow \infty$ and using that $\tilde{\xi}(x_1)$ decays as $|x_1|^{\eta + \epsilon ' - 1}$, we get that
for any $b$ that
$$ \sup_K I_K \geq C'' \int_b^a x_1^{\eta + \epsilon ' - \epsilon_1 - 1}|\ln x_1|^d  \,dx_1 \eqno (3.13)$$
Since $ \sup_K I_K$ is finite, we must therefore have that $\eta + \epsilon ' - \epsilon_1 > 0$, contradicting the choice of 
$\eta$. Hence we have arrived at a contradiction and the proof of part b) of Theorem 2.2 is complete, thereby completing the
proof of the whole theorem.

\noindent In the proof of Theorem 2.3 we will use the following lemma.

\noindent {\bf Lemma 3.1.} Given any multiindex $(a,b)$ there is a neighborhood $N$  of the origin and a constant $C_{a,b,f,N}$
such that on $N$ one has 
$$|\partial_{x_1}^a\partial_{x_2}^b f(x_1,x_2)| \leq C_{a,b,f,N}{1 \over |x_1|^a|x_2|^b} f^*(x_1,x_2) \eqno (3.14)$$
If each $f_e(x_1,x_2)$ has no zeroes on $(\R - \{0\})^2$, then  there is in addition a neighborhood $N'$ 
of the origin and a constant $c_{f,N'}$ such that $|f(x_1,x_2)| \geq c_{f,N'}f^*(x_1,x_2)$ on $N'$. In other words,
$|f(x_1,x_2)|  \sim f^*(x_1,x_2)$ on $N'$.

\noindent {\bf Proof.} As mentioned in the proof of Lemma 1.2,  Lemma 2.1 of [G2] implies that for any real analytic function
 $g(x_1,x_2)$  on a neighborhood of the origin  with $g(0,0) = 0$, there is an inequality 
$|g(x_1,x_2)| \leq Cg^*(x_1,x_2)$ on a neighborhood of the origin.
Applying this to any $\partial_{x_1}^a\partial_{x_2}^b f$ equal to zero at origin gives $(3.14)$ for that $(a,b)$. If $\partial_{x_1}^a\partial_{x_2}^b f(0,0) \neq 0$ the inequality is immediate, so $(3.14)$ holds in all cases. 

Suppose now that each $f_e(x_1,x_2)$ has no zeroes on $(\R - \{0\})^2$. We divide a neighborhood $N$ of the origin into wedges
$A_i$ and $B_i$ as follows. Each $A_i$ is of the form $\{ (x_1,x_2) \in N: {1 \over K}|x_1|^{m_i} < |x_2| < K|x_1|^{m_i}\}$ for some
large $K$ and where $e_i$ is a compact edge of $N(f)$. Each $B_i$ is of the form  $\{ (x_1,x_2) \in N:  K|x_1|^{m_i} < |x_2| 
< {1 \over K}|x_1|^{m_{i+1}}\}$ for compact edges $e_i$ and $e_{i+1}$ of $N(f)$,  or is of the form
$\{ (x_1,x_2) \in N:  K|x_1|^{m_{n-1}} < |x_2| \}$, or is of the form $\{ (x_1,x_2) \in N: |x_2| < {1 \over K}|x_1|^{m_1}\}$.

In the setting of Lemma 2.1 of [G2], the $A_i$ and $B_i$ are the sets denoted by $W_{ij}$. For the $A_i$, Lemma 2.1 of [G2]
says that  given any fixed $K$ and any $\delta > 0$, there is a neighborhood $V_i$ of the origin such that on $A_i \cap V_i$ we have
$$|f(x_1,x_2) - f_e(x_1,x_2)| < \delta |x_1|^{v_1^i}|x_2|^{v_2^i} \eqno (3.15)$$
 In addition, using that  that each $f_e(x_1,x_2)$ has no zeroes on $(\R - \{0\})^2$, the mixed homogeneity of $f_e(x_1,x_2)$,
and the resulting fact that $f_e(1,x)$ and $f_e(1,-x)$ 
have no zeroes on $[{1 \over K}, K] \cup [-K, -{1 \over K}]$, there is a constant
$c$ such that $|f_e(x_1,x_2)| > c |x_1|^{v_1^i}|x_2|^{v_2^i}$ on $A_i$. Hence choosing $\delta = {c \over 2}$, we 
conclude that
$|f(x_1,x_2)| > {c \over 2} |x_1|^{v_1^i}|x_2|^{v_2^i}$ on $A_i \cap V_i$. By
$(2.4)$ we have that  $|x_1|^{v_1^i}|x_2|^{v_2^i} \sim f^*(x_1,x_2)$ on $A_i \cap V_i$, so there exists a  $c'$ for which 
$|f(x_1,x_2)| > c' f^*(x_1,x_2)$ on $A_i \cap V_i$ as needed.

For the $B_i$, Lemma 2.1 of [G2] says that there is a single vertex $v_j$ of $N(f)$ and a $d_j \neq  0$ such that given any $\delta > 0$,  if $K$ were chosen large enough  there is a neighborhood $U_i$ of the origin such that
$|f(x_1,x_2) -  d_j v_1^j v_2^j| < \delta |x_1|^{v_1^j}|x_2|^{v_2^j}$ on $B_i \cap U_i$.
Taking $\delta < {1 \over 2}|d_j|$, we have $|f(x_1,x_2)| >  {1 \over 2}|d_j v_1^j v_2^j|$ on $B_i \cap U_i$. By $(2.4)-(2.6)$ we have that $|v_1^j v_2^j| \sim f^*(x_1,x_2)$ on $B_i$, so $|f(x_1,x_2)| >  c''f^*(x_1,x_2)$ on
$B_i \cap U_i$ for some constant $c'' > 0$ as needed. 

Letting $N'$ be the intersection of all $U_i$ and $V_i$, we see that
on $N'$ we have $|f(x_1,x_2)| > c_{f,N'}f^*(x_1,x_2)$ for some constant $c_{f,N'}$ as needed. This completes the proof of
 Lemma 3.1.

\noindent {\bf Proof of Theorem 2.3.}

\noindent We write $F(\lambda_1,\lambda_2) = \sum_{j,k}F_{jk}(\lambda_1, \lambda_2)$, where 
$$F_{jk}(\lambda_1, \lambda_2) = \int_{\R^2} \phi(x_1,x_2)\beta(2^j x_1)\beta(2^k x_2) 
|f(x_1,x_2)|^{-\rho}e^{-i\lambda_1x_1 - i\lambda_2x_2}\,dx_1\,dx_2 \eqno (3.16)$$
Here $\beta(x)$ is a nonnegative smooth compactly supported function on $\R$ whose support  does not intersect
 some neighborhood of $0$.
If $j$ is such that $2^{-j} > |\lambda_1|^{-1}$, we integrate by parts in $(3.16)$ in the $x_1$ variable $N$ times, integrating
the $e^{-i\lambda_1x_1}$ and differentiating the rest. Each time we do so we get a ${1 \over i\lambda_1}$ from the integration.
When the derivative lands on $\beta(2^j x_1)$ or one of its derivatives we get a factor of $C2^j$, and each time the 
derivative lands on  $\phi(x_1,x_2)\beta(2^k x_2)$ or one of its derivatives we also get a factor of $C2^j$ due to the conditions
$(2.8)$.

As for when the derivative lands on the
$ |f(x_1,x_2)|^{-\rho}$ factor, by Lemma 3.1 $f(x_1,x_2)$ is of a single sign in each quadrant, so $|f(x_1,x_2)|^{-\rho}$ is
$(\pm f(x_1,x_2))^{-\rho}$ on a given quadrant. Each time the $x_1$ derivative lands on such a factor or one of its derivatives,
by $(3.14)$ and the fact that $f(x_1,x_2) \sim f^*(x_1,x_2)$, 
one gets a factor bounded by $C{1 \over |x_1|}$.  Due to the support conditions on $\beta(2^j x_1)$, this too is bounded by
$C2^j$. 

We conclude that the integration by parts leads to an overall factor of $C2^j$. Hence $N$ integrations by parts results in a factor of $C2^{jN} \sim C|x_1|^{-N}$. Because of this and the fact that
$f(x_1,x_2) \sim f^*(x_1,x_2)$ by Lemma 3.1, we conclude that there is a neighborhood $U$ of the origin such that if 
$\phi(x_1,x_2)$ is supported in $U$, then for any $N$ we have an estimate
$$|F_{jk}(\lambda_1, \lambda_2)| \leq C_N {1 \over |\lambda_1|^N}\int_{U \cap \{x: 2^{-j-1} < |x_1| < 2^{-j},\, 2^{-k -1} < |x_2| < 2^{-k}\}}
{1 \over |x_1|^N}|f^*(x_1,x_2)|^{-\rho}\,dx_1\,dx_2  \eqno (3.17a)$$
In exactly the same way, reversing the roles of the $x_1$ and $x_2$ variables, if $2^{-k} > |\lambda_2|^{-1}$ we have
$$|F_{jk}(\lambda_1, \lambda_2)| \leq C_N {1 \over |\lambda_2|^N}\int_{U \cap \{x: 2^{-j-1} < |x_1| < 2^{-j},\, 2^{-k -1} < |x_2| < 2^{-k}\}}
{1 \over |x_2|^N}|f^*(x_1,x_2)|^{-\rho}\,dx_1\,dx_2  \eqno (3.17b)$$
If both $2^{-j} > |\lambda_1|^{-1}$ and $2^{-k} > |\lambda_2|^{-1}$, we can first do $N$ integrations by parts in the $x_1$
variable followed by $N$ integrations by parts in the $x_2$ variable to obtain that $|F_{jk}(\lambda_1, \lambda_2)|$ is bounded by
$$C_N {1 \over |\lambda_1\lambda_2|^N}\int_{U \cap \{x: 2^{-j-1} < |x_1| < 2^{-j},\, 2^{-k -1} < |x_2| < 2^{-k}\}}
{1 \over |x_1 x_2|^N}|f^*(x_1,x_2)|^{-\rho}\,dx_1\,dx_2  \eqno (3.17c)$$
We will obtain our desired estimates for a given $F_{jk}(\lambda_1,\lambda_2)$ as follows. When $2^{-j} \leq |\lambda_1|^{-1}$
and $2^{-k} \leq  |\lambda_2|^{-1}$ we just take absolute values in $(3.16)$ and integrate. When $2^{-j} > |\lambda_1|^{-1}$
and $2^{-k} < |\lambda_2|^{-1}$ we use $(3.17a)$. When $2^{-j} < |\lambda_1|^{-1}$
and $2^{-k} > |\lambda_2|^{-1}$ we use $(3.17b)$, and when $2^{-j} \geq |\lambda_1|^{-1}$
and $2^{-k} \geq |\lambda_2|^{-1}$ we use $(3.17c)$. We will add over all $j$ and $k$ to obtain the desired estimates.
The value of $N$ will be determined  by our arguments. 

Taking absolute values in $(3.16)$, integrating, and adding over all $j$ and $k$ with $2^{-j} \leq |\lambda_1|^{-1}$
and $2^{-k} \leq |\lambda_2|^{-1}$ leads to the desired estimate
$$C\int_{\{x \in U:\, |x_1| < |\lambda_1|^{-1}, \,|x_2| < |\lambda_2|^{-1}\}}|f^*(x_1,x_2)|^{-\rho}\,dx_1\,dx_2  \eqno (3.18)$$
We next add over all $(j,k)$ such that $2^{-j} > |\lambda_1|^{-1}$
and $2^{-k} < |\lambda_2|^{-1}$. For a given $k$, we add
 estimates $(3.17a)$ in $j$. Let $a$ denote the minimum $v_1^i$ appearing in any of the terms $x_1^{v_1^i}x_2^{v_2^i}$ 
defining $f^*(x_1,x_2)$. Then $f^*(2x_1,x_2) \geq 2^a f^*(2x_1,x_2)$, and $(f^*(2x_1,x_2))^{-\rho} \leq 
2^{-\rho a} (f^*(2x_1,x_2))^{-\rho}$. (If $\rho < 0$, we let $a$ be the maximal $v_1^i$). Thus if $N$ is large enough, the 
integrand in $(3.17a)$ decreases by a factor of at least $4$ each time $j$ increases by $1$ for fixed $k$. As a result, the integral
decreases by a factor of at least $2$ each time $j$ increases by $1$ for fixed $k$. Hence the sum of $(3.17a)$ over $j$ with 
$2^{-j} > |\lambda_1|^{-1}$ is bounded by a constant times what one gets in $(3.17a)$ setting $|\lambda_1| =  2^{-j}$, namely
$$C\int_{\{x \in U:\, {1 \over 2}|\lambda_1|^{-1} < |x_1| < |\lambda_1|^{-1}, \,2^{-k -1} < |x_2| < 2^{-k}\}}|f^*(x_1,x_2)|^{-\rho}\,dx_1\,dx_2  \eqno (3.19)$$
We now add $(3.18)$ over all $k$ with  $2^{-k} < |\lambda_2|^{-1}$, and we see that the sum of all terms with 
$2^{-j} > |\lambda_1|^{-1}$ and $2^{-k} < |\lambda_2|^{-1}$ is bounded by
$$C\int_{\{x \in U:\,{1 \over 2}|\lambda_1|^{-1} < |x_1| < |\lambda_1|^{-1}, \,|x_2| < |\lambda_2|^{-1}\}}|f^*(x_1,x_2)|^{-\rho}\,dx_1\,dx_2  \eqno (3.20)$$
This too is bounded by the desired estimate $(3.18)$, so we are done with the terms where $2^{-j} > |\lambda_1|^{-1}$
and $2^{-k} < |\lambda_2|^{-1}$. By symmetry, the same method gives this estimate for the terms where $2^{-j} <
 |\lambda_1|^{-1}$ and $2^{-k} > |\lambda_2|^{-1}$, using $(3.17b)$ in place of $(3.17a)$. 

It remains to consider the terms where $2^{-j} \geq |\lambda_1|^{-1}$ and $2^{-k}  \geq |\lambda_2|^{-1}$.
 If $N$ is large enough,
similar to the argument leading to $(3.19)$, increasing $j$ by 1 in the expression $(3.17c)$ for a fixed $k$ decreases the term 
by a factor of at least 2. Hence adding over all such $j$ leads to a bound of a constant times the term where $j$ is minimal, that is
where $2^{-j}$ is within a factor of 2 of $|\lambda_1|^{-1}$. But this term is exactly $(3.17b)$ with $2^{-j} = |\lambda_1|^{-1}$. Hence adding these 
over all $k$ with $2^{-k}  \geq |\lambda_2|^{-1}$ once again leads to the desired bound $(3.18)$. This completes the proof of 
Theorem 2.3.

\noindent {\bf References.}

\noindent [AGV] V. Arnold, S. Gusein-Zade, A. Varchenko, {\it Singularities of differentiable maps},
Volume II, Birkhauser, Basel, 1988. \parskip = 4pt\baselineskip = 4pt

\noindent [G1] M. Greenblatt {\it  Newton polygons and local integrability of negative powers of
smooth functions in the plane}, Trans. Amer. Math. Soc., {\bf 358} (2006), 657-670.

\noindent [G2] M. Greenblatt, {\it Oscillatory integral decay, sublevel set growth, and the Newton
polyhedron}, Math. Annalen {\bf 346} (2010), no. 4, 857-895.

\noindent [G3] M. Greenblatt,  {\it Singular integral operators with kernels associated to negative powers of real analytic functions},
 J. Funct. Anal. {\bf 269} (2015), no. 11, 3663-3687. 

\noindent [G4] M. Greenblatt, {\it Convolution kernels of 2D Fourier multipliers based on real analytic functions}, preprint.

\noindent [G5] M. Greenblatt, {\it Uniform bounds for Fourier transforms of surface measures in $\R^3$ with nonsmooth density}, to appear, Trans. AMS.

\noindent [V] A. N. Varchenko, {\it Newton polyhedra and estimates of oscillatory integrals}, Functional 
Anal. Appl. {\bf 18} (1976), no. 3, 175-196.

\end